\documentclass[a4paper,11pt,leqno]{article} 

\usepackage{lineno}

\usepackage{amsmath,amsfonts} 
\usepackage{graphics}
\usepackage[T1]{fontenc}
\usepackage{lineno}
\usepackage[utf8]{inputenc}
\usepackage{amsmath,mathrsfs,amssymb,amsfonts,nicefrac}
\usepackage[english]{babel}
\usepackage[active]{srcltx}
\usepackage{tikz}
\usepackage{dsfont}
\usetikzlibrary{patterns}

\numberwithin{equation}{section}

\newtheorem{thm}{Theorem}[section]

\newtheorem{lm}[thm]{Lemma}  

\newtheorem{rem}[thm]{Remark}  

\newcommand{\RR}{\mathbb{R}}   
   
\newcommand{\di}{\displaystyle}

\renewcommand{\epsilon}{\varepsilon}
\newcommand{\e}{\varepsilon}

\begin{document}   
   
\title{\textbf{Interaction of a  free boundary with a diffusion on a plane: analogy with the obstacle problem}}

\author{
{\bf Luis A. Caffarelli} \\ 
The University of Texas at Austin\\
Mathematics Department RLM 8.100, 2515 Speedway Stop C1200\\
 Austin, Texas 78712-1202, U.S.A.  \\
\texttt{caffarel@math.utexas.edu}
\\[2mm]
{\bf Jean-Michel Roquejoffre} \\ 
Institut de Math\'ematiques de Toulouse (UMR CNRS 5219) \\ 
Universit\'e Toulouse III-Paul Sabatier,
118 route de Narbonne\\
31062 Toulouse cedex, France \\ 
\texttt{jean-michel.roquejoffre@math.univ-toulouse.fr}
\\[2mm]
{\bf Ignacio Tomasetti}\\
The University of Texas at Austin\\
Mathematics Department RLM 8.100, 2515 Speedway Stop C1200\\
 Austin, Texas 78712-1202, U.S.A.  \\
\texttt{tomasetti@math.utexas.edu}
}
  \date{}
\maketitle   
\begin{abstract} 
\noindent  A horizontal $N$-dimensional plane, having a diffusion
of its own, exchanges with the lower half space. There, a reaction-diffusion
process, modelled by a free boundary problem, takes place. We wish to
understand whether, and how, the free boundary meets the plane.

\noindent The origin of this problem is a two-dimensional reaction diffusion model proposed
some time ago by the second author, in collaboration with H. Berestycki and L. Rossi,
 to model how biological
invasions can be enhanced by a line of fast diffusion. Some counter-intuitive numerical simulations of this model, 
due
to A.-C. Coulon,  have been explained by the first two authors   by transforming the model into a free boundary
interacting with a line, and a careful study of the free boundary. At this occasion,
it was noticed that the free boundary very much like that of the obstacle problem.

\noindent The goal of the paper is to explain how this analogy with the obstacle problem can
be pushed further in higher space dimensions.

\end{abstract}

\section{Introduction}    

\subsection{Model and question}
Let $\varphi$ be a   smooth concave function of $\RR^N$, attaining a positive maximum at $x=0$.
We look for a function $u(X)$, defined for $X=(x,y)\in\RR^{N+1}_-$ with $\RR^{N+1}_-:=\RR^N\times\RR_-$, and a surface $\Gamma\subset\RR^{N+1}_-$ such that
\begin{equation}
\label{e1.1}
\left\{
\begin{array}{rll}
\Delta u=&0\quad(x,y)\in\{u>0\}\\
\vert\nabla u\vert=&1\quad((x,y)\in\Gamma:=\partial\{u>0\})\\
\ \\
-\Delta_x u+ u_y=&0\quad(x\in\RR^N, y=0)\\
u(x,0)\geq&\varphi(x).
\end{array}
\right.
\end{equation}
We ask whether, and how, the free surface $\Gamma$ meets the plane $\{y=0\}$.
\subsection{Motivation of the model, possible extensions}
The study of Model \eqref{e1.1} is motivated by previous investigations on a two-dimensional version that modelled the influence of a line of fast diffusion on the propagation of reaction-diffusion waves. More precisely, in \cite{CR1} and \cite{CR2},
the first two authors study o a system of the form
 \begin{equation}
\label{e1.110}
\left\{
\begin{array}{rll}
-\Delta u+c\partial_x u=&0\quad(x,y)\in\{u>0\}\\
\vert\nabla u\vert=&1\quad((x,y)\in\Gamma:=\partial\{u>0\}\\
\ \\
-u_{xx}+c\partial_xu+ u_y=&0\quad\hbox{for $x\in\RR$, $y=0$}\\
u(-\infty,y)=&1,\quad\hbox{uniformly in $y\in\RR_-$}\\
 u(+\infty,y)=&0\quad\hbox{pointwise in $y\in\RR_-$}.
\end{array}
\right.
\end{equation}
Studying \eqref{e1.110} was itself motivated by a rather a counter-intuitive numerical simulation in \cite{ACC}, of a system initially proposed in \cite{BRR} to model the influence of transportation networks on biological invasions. See also \cite{BCRR}. 

\noindent Among other things, it is proved in \cite{CR1} that the free boundary $\Gamma$ hits the line $\{y=0\}$ as an inverted parabola: assuming that $(0,0)$ is the intersection point, then, in its vicinity, $\Gamma$ may be described as
$$
y=-\frac{x^2}2,\ \ \ \ x\leq0,
$$
and $u(x,0)$ grows in a quadratic fashion for negative $x$. This is typical of the behaviour of the solutions to the obstacle problem, and the present study of Model \eqref{e1.1} is an attempt to understand whether the situation discovered in \cite{CR1} is general.
Theorems \ref{t1.1} and \ref{t1.2} below show that this is indeed the case.

\noindent We point out that, while Model shows the genericity of the two-dimensional situation in several space dimensions, our study depends on the global geometry of the problem, and that properties such as monotonicity in every direction, help simplifying 
the analysis. Thus, a general local study of the interaction of the free boundary of the one phase problem with the plane is still to be carried out. Whether, in particular, Theorem \ref{t1.2} would persist in a general setting, is something we do not know. Another interesting question is to replace the equation for $u$ on the plane by  the family of equations 
$$
-D\Delta u+u_y=0,
$$
and investigate what happens as $D\o+\infty$. This question, by the way, is studied in the forthcoming work \cite{CR3}.

\noindent Other geometries would be relevant. In particular, the obstacle condition
$$
u(x,0)\geq\varphi(x)
$$	
could certainly be replaced by a Dirichlet type condition on a vertical cylinder inside the lower half space. The ideas displayed here would probably work in this setting without so many changes. The diffusion equation on the plane $\{y=0\}$ could also be replaced by
$$
-\sum_{1\leq i,j\leq N}\partial_{x_i}\bigl(a_{ij}(\nabla u)\bigl)+u_y=0,
$$
where the matrix $\bigl(a_{ij}(p)\bigl)$ satisfies the usual ellipticity conditions, such as, for instance, in  \cite{BK}. This last extension could open the door to the study of interesting directional effects.
\subsection{Main results}
\noindent Let us first introduce some notations.  As we will have to deal with several choices of balls, either in the lower half space $\RR^{N+1}_-$ or the plane $\{y=0\}$ it is convenient to denote them with distinctive signs now. So, let $B^{N+1}_R(X)$ be the ball of $\RR^{N+1}$ with centre $X$ and radius $R$; when $X=0$ we will simply denote it by $B^{N+1}_R$. In the same spirit, we denote by $B_R^N(x)$ the ball of the plane $\{y=0\}$ with centre $x$  and radius $R$, again denoted by $B_R^N$ if $x=0$.

\noindent  Let us now make a set of slightly more stringent assumptions on $\varphi$. They are probably not indispensable, however they will ensure, in a relatively painless way, the existence of a solution to \eqref{e1.1} whose support on the plane is larger than that of $\varphi$.  
\begin{enumerate}
\item We have $\varphi(x)=1$ if $x\in B_1^N$.
\item There is $\rho_0>0$ such that
\begin{equation}
\label{e1.3}
\mathrm{supp}~\varphi^+\subset B^N_{1+\rho_0}.
\end{equation}
Moreover $\varphi$ is strictly concave outside $B^N_1$.
\end{enumerate}
With these assumptions we can prove the 
\begin{thm}\label{t1.1}
If $\rho_0>0$ is small enough, Problem \eqref{e1.1} has a solution $(u,\Gamma)$ which has the following properties: 
\begin{enumerate}
\item there is $\rho_1>\rho_0$ such that  
\begin{equation}
\label{e1.4}
B^N_{1+\rho_1}\subset\mathrm{supp}~u(.,0).
\end{equation}
\item There is $\theta_0>0$ such that, for every $x\in\RR^N$ outside $\mathrm{supp}~\varphi^+$,  we have $\partial_eu(x,0)\leq0$ for every direction $e$ making an angle less or equal $\theta_0$ with $\di\frac{x}{\vert x\vert}$. 
\item The free boundary $\Gamma$ meets the plane $\{y=0\}$ along a Lipschitz surface, and encloses a domain $\Omega$ of the plane $\{y=0\}$.
\end{enumerate}
\end{thm}
This result is  in the spirit of the work \cite{C00} of the first author, where some a priori monotonicity results of the solution entail the Lipschitz property of the free boundary in an obstacle type problem.
\begin{thm}
\label{t1.2}
There is $\delta>0$ such that  have 
$$
u_y(x,0)\geq\delta  \quad\hbox{if $x\in\Omega$}.
$$
Moreover,  $u_y$ is H\"older-continuous in $\bar\Omega$.
\end{thm}
The end result of Theorem \ref{t1.2} is that   $u(x,0)$ behaves, in the vicinity of $\Gamma\cap\{y=0\}$, very much like the solution of the obstacle problem in its reduced form, that is, according to the notation of \cite{C1} (Part I, Point b)):
$$
-\Delta_xu(x,0)+g(x)H(u)=0\quad\hbox{in $\Omega=\{u(x,0)>0\}$},
$$
where $H$ is the Heaviside function, and the function $g(x)$ is a positive, H\"older continuous function, a role that the function $u_y(x,0)$, can be made to play   thanks to Theorem \ref{t1.2}. This property is sufficient to study a good part of the regularity of $\partial\Omega$, see for instance  \cite{C1}, or the original paper \cite{C0}, or the review \cite{C0review}. In particular,  Theorem 6 of \cite{C1} implies that $\partial\Omega$ is a $C^{1,\alpha}$ surface. 
\subsection{Organisation of the paper}
A first issue is to construct a solution to \eqref{e1.1} having some nice geometrical properties, such as the monotonicity in a cone of directions of the plane. While this should not be indispensable to the study (at least the basic properties), it greatly helps us in  the analysis. And so, the next two sections are devoted to the construction of such a solution to \eqref{e1.1}, through an approximation by a singularly perturbed semilinear problem. We prove, in particular, that the free boundary $\Gamma$ intersects the plane, opening the way to the analogy with the obstacle problem. In Section 4, we start the analysis of $u(.,0)$ by proving the basic properties of the solutions to the obstacle problem. We prove optimal regularity, that is, $u(.,0)$ is $C^{1,1}$ in the vicinity of 
$\Gamma\cup\{y=0\}$, and that $u(x,.)$ grows at least quadratically as $x$ departs from the free boundary. The last section is devoted to the H\"older continuity of $u$, thus concluding the analogy.
\section{Construction of a solution to \eqref{e1.1}: the semilinear approximation}
The goal of the next two sections is to construct a solution $(u,\Gamma)$ of \eqref{e1.1}  that is as classical as possible, that is  the free boundary relation is satisfied in the strong sense away from the place where $\Gamma$ meets the plane, and the exchange condition $-\Delta_xu+u_y=0$ is satisfied in the strong sense, the function $u_y$ being $L^\infty$ at this stage. For that we will resort to a smoothing of the free boundary problem by a semilinear one. This procedure, currently used in mathematical models for premixed flames, was first used in multidimensional elliptic problems by Berestycki, Nirenberg and the first author in \cite{BCN} to approximate one-phase free boundary problems. See also the work of Lee, Mellet and the first author where this approximation is used in the homogenisation of this type of free boundary problems for parabolic equations. The procedure  was more recently used in \cite{CR1} to construct a travelling wave solution.

\noindent As is usual, the coincidence set of $u$ will be the (closed) subset of the plane $\RR^N$ where $u(.,0)=\varphi$. We denote it by $\mathcal{C}_u$.

\noindent Consider a smooth function  $f$, positive on $\RR_+$, and such that $\di\int\beta(u)du=\di\frac1{\sqrt2}$, with $\beta(u)=uf(u)$, and such that $uf'(u)+f(u)\geq0$ for $u\in\RR_+$. A suitable multiple of $\di\frac1{1+u}$ will do the trick. 

\noindent We will first solve the problem
\begin{equation}
\label{e2.1}
\left\{
\begin{array}{rll}
-\Delta u+\beta(u)=&0\quad(x,y)\in\RR^{N+1}_-\\
\ \\
-\Delta_x u+ u_y=&0\quad(x\in\RR^N, y=0)\\
u(x,0)\geq&\varphi(x).
\end{array}
\right.
\end{equation}

\noindent Once a solution  $u$ is constructed, we will set, for every $\e>0$:
\begin{equation}
\label{e2.8}
\beta_\e(u)=\di\frac{u}{\e^2}f(\frac{u}\e),
\end{equation}
and show that, up to a subsequence, the family $(u_\e)_\e$ will converge, as $\e>0$, to a solution of \eqref{e1.1}.

\noindent Let us therefore start with \eqref{e2.1}.  We start the construction in a standard way: we first look for a solution to the approximate problem
\begin{equation}
\label{e2.3}
\left\{
\begin{array}{rll}
-\Delta u+\beta(u)=&0\quad(x,y)\in\RR^{N+1}_-\cap B^{N+1}_R\\
-\Delta_x u+ u_y=&0\quad(x\in\RR^N, y=0)\cap B^N_R\\
\ \\
u(x,0)\geq&\varphi(x)\\
u(X)=&0\quad(\vert X\vert=R).
\end{array}
\right.
\end{equation}
In what follows, $R$ is assumed to be large enough so that $\mathrm{supp}~\varphi^+$  is inside $B_{R/2}^N$ We have the  easy a priori estimate
$$
0\leq u\leq \varphi(0).
$$
The energy functional
\begin{equation}
\label{e2.9}
J_R(u)=\int_{B_R^{N+1}\cap\RR^{N+1}_-}\biggl(\vert\nabla u\vert^2+B(u)\biggl)dX+\int_{B_R^N}\vert\nabla_xu(x,0)\vert^2dx,
\end{equation}
with
$$
\quad B(u)=\int_0^u\beta(v)dv,
$$
is easily shown to attain a minimum over all functions $u\in H^1_0(B^{N+1}_R\cap\RR^{N+1}_-)$, whose trace on the plane $\{y=0\}$ is larger than $\varphi$ and belongs to $H^1_0(B^N_R)$. 
This means that (\ref{e2.3}) has at least a solution. Now, we need to prove enough smoothness - not such a difficult issue at this stage. The main step is a uniform Lipschitz bound.
\begin{lm}
\label{l2.1}
There is a constant $M>0$, independent of $\lambda$ and $R$  (but, at this stage, strongly depending on $\beta$) such that, for all $u\in S_{\lambda,R}$ we have
$$
\Vert\nabla u\Vert_{L^\infty(B_R^{N+1}\cap\RR^{N+1}_-)}\leq M.
$$
\end{lm}
\noindent {\sc Proof.} The only point that is slightly nonclassical here is to handle the detachment from the coincidence set $\mathcal{C}_u$, and the points where the Wentzell conditions are enforced at the same time. In order to isolate these two parts, let $R_0>0$ be such that $B_{R_0}^N$ contains the whole coincidence set, and choose $R\geq 3R_0$ once and for all.
Outside $B_{3R_0}^{N+1}\cap\RR^{N+1}_-$, the problem is a standard semilinear equation with Wentzell boundary conditions, so that $u$ enjoys $C^{2,\alpha}$ bounds that are uniform in $R$.

\noindent To deal with the remaining part, we could use the property $uf'(u)+f(u)\geq0$ for a quick proof. The proof that we are going to give, while being slightly longer, shows that this particular property is not needed here.  Let us first recall some basic facts about the first eigenvalue. If $\mathcal{O}$ is a bounded open subset of $\RR^{N+1}_-$ with trace $\mathcal{O}_0$ on the plane $\RR^N$, contained in $B_{2R_0}$, and such that $\mathcal{C}_u$  is contained in $\mathcal{O}_0$, at nonzero distance from $\partial\mathcal{O}_0$. Consider the eigenvalue problem
\begin{equation}
\label{e2.6}
\left\{
\begin{array}{rll}
-\Delta u=&\mu u\quad((x,y)\in\mathcal{O})\\
-\Delta_xu+u_y=&\mu u\quad(x\in\mathcal{O}_0,y=0)\\
\ \\
u(x,y)=&0\quad((x,y)\in\mathcal{O},y<0)\\
u(x,0)=&0\quad(x\in\partial\mathcal{O}_0).
\end{array}
\right.
\end{equation}
 From Krein-Rutman's theorem there is a smallest eigenvalue, called $\mu_1(\mathcal{O})$, associated to a positive eigenfunction $\Phi_1(x,y)$, and characterised as the minimum of the Rayleigh type quotient
$$
\frac{\di{\int_{\mathcal{O}_0}\vert\nabla u(x,0)\vert^2dx+\int_{\mathcal{O}}\biggl(\vert\nabla u\vert^2+u^2\biggl)dxdy}}{\Vert u\Vert_{L^2(\mathcal{O}_0)}^2+\Vert u
\Vert_{L^2(\mathcal{O})}^2},
$$
taken on all functions $u\in H^1(\mathcal{O})$,whose trace on $\partial\mathcal{O}\cap\{y<0\}$ is zero,  and such that $u(.,0)\in H^1_0(\mathcal{O}_0)$. Let us define the width of $\mathcal{O}$ as the smallest $l>0$ such that there is a
 point $(x,-l)$ in $\mathcal{O}$, one easily deduces from the Poincar\'e inequality that $\mu_1(\mathcal{O})$ goes to infinity as the width of $\mathcal{O}$ goes to 0.

\noindent Choose $\mathcal{O}$ and $l_0>0$ such that 
\begin{itemize}
\item $\mathcal{O}$ is a smooth subdomain of $B^{N+1}_{2R_0}\cap(-l_0,0)$, 
\item $\partial\mathcal{O}$ contains $B_{4R_0/5}^N\times\{-l_0\}$ and $\mathcal{O}_0$ contains $B_{4R_0/5}^N$,
\item  we have $\mu_1(\mathcal{O})\geq\Vert\partial_u\beta(u)\Vert_\infty+1$.
\end{itemize}
Choose, finally, a subdomain $\mathcal{O}'$ of $\mathcal{O}$ such that $\partial\mathcal{O}'\cap\{y<0\}$  is 
at positive distance from $\partial\mathcal{O}\cap\{y<0\}$, and such that, if $\mathcal{O}_0'$  is the trace of $\mathcal{O}_0$ on the plane $\RR^N$, then 
$\partial\mathcal{O}_0'$ is at a positive distance from $\partial\mathcal{O}_0$.

\noindent Pick any direction $e\in\RR^N$, let us bound $\partial_eu$. As mentionned before, it is sufficient to do it in $\mathcal{O}'$.  If $\Phi_1$ is the eigenfunction associated to
$\mu_1(\mathcal{O})$, then $C\Phi_1$ bounds $\vert\partial_eu\vert$ on $\mathcal{O}'\cap\{y<0\}$, as well as on $\mathcal{C}_u$, as soon as $C>0$ is large enough, in any case larger than $\Vert\nabla\phi\Vert_\infty$. Now, $v=\partial_eu$ satisfies
$$
-\Delta v+\partial_u\beta(u)v=0.
$$
Since $\mu_1(\mathcal{O}')>\mu_1(\mathcal{O})\geq\Vert\partial_u\beta(u)\Vert_\infty+1$, the maximum principle holds and we have $\vert v\vert\leq C\vert\Phi_1\vert$. The same
argument applies to $w:=\partial_{ee'}u$, where $e'$ is another direction of $\RR^N$. We have indeed
$$
-\Delta w+\partial_u\beta(u)w=\partial_{uu}\beta(u)\partial_eu\partial_{e'}u=O(1).
$$
It remains to bound $u_y$. We denote it by $v$, and notice that $v$ is bounded on the plane outside $\mathcal{C}_u$, while $v_y$ is bounded on $\mathcal{C}_u$, as it is equal to $-\Delta\phi$. Consider $v(x,y)+y\Delta\varphi(x)$ and note that this time it vanishes on $\mathcal{C}_u$, while still being bounded outside. The preceding argument may be applied to $v$, this time to the Laplacian in $\mathcal{O}'$, still with Dirichlet conditions on $\mathcal{O}'\cap\{y<0\}$, but with Neumann conditions on $\mathcal{O}_0'$.
\hfill$\Box$

\noindent As a consequence of the equation $-\Delta u+u_y=0$, for all $\alpha\in(0,1)$, there is a uniform $C^{1,\alpha}$ on $B^N_R$   for $u(.,0)$, for all $u\in S_{\lambda,R}$. This bound is in turn translated into a uniform 
$C^{1,\alpha}$ bound in $B_R^{N+1}\cap\RR^{N+1}_-$.  Once this is at hand, we go to the main qualitative property. 
\begin{lm}
\label{l2.2}
There is $\delta>0$, independent of $\lambda$ and $R$, such that we have, for  all solution $u$ of (\ref{e2.3}),
for all $x$ outside $\mathrm{supp}~\varphi^+$, and all unit vector $e$:
\begin{equation}
\label{e2.20}
\hbox{If $\di{\biggl\vert\frac{x}{\vert x\vert}-e\biggl\vert}\leq\delta$, then}\ \partial_eu(x,0)\leq0.
\end{equation}
Moreover we have
\begin{equation}
\label{e2.5}
 \partial_yu\geq\delta.
\end{equation}
\end{lm}
\noindent {\sc Proof.}   This is really where we need $uf'(u)+f(u)\geq0$. Let us first check (\ref{e2.5}) in the most standard way, by proving that $v:=u_y$ cannot have a nonpositive minimum. Compactness will once again ensure $v\leq -\delta$, for some  $\delta>0$ independent of $R$. Assume therefore the existence of a nonpositive minimum for $v$. As $uf'(u)+f(u)\geq0$, the  equation for $v$
$$
-\Delta v+\bigl(uf'(u)+f(u)\bigl)v=0\quad\hbox{in $B_R^{N+1}\cap\RR^{N+1}_-$}
$$
and the Hopf Lemma ensures that the minimum can be assumed neither inside $B_R^{N+1}\cap\RR^{N+1}_-$, nor on $\partial(B_R^{N+1}\cap\RR^{N+1}_-)$. So, it has to be taken on the plane. However there are new obstructions: on $\mathcal{C}_u$  we have
$$
v_y(x,0)=u_{yy}(x,0)=-\Delta_xu(x,0)=-\Delta\varphi(x)>0,
$$
so that the Hopf Lemma precludes a minimum of $v$ on that set. Outside we have, by the same considerations:
$$
v_y(x,0)+v(x,0)=\beta\bigl(u(x,0)\bigl),
$$ 
which forbids a nonpositive minimum once again.

\noindent That $u$ decreases along each direction $e$  of the plane $\{y=0\}$, outside $\mathrm{supp}~\varphi^+$, is a standard consequence of the moving plane method \cite{GNN}, as one can perform it by reflecting $u$ across any hyperplane 
\begin{equation}
\label{e2.15}
H_\lambda=\{(x,y):\ x.e=\lambda\}.
\end{equation}
 Indeed, if $u^r$ is the reflection of $u$, we have $-\Delta u^r+\beta(u^r)=0$ in the reflection of $B_R^{N+1}\cap\{y<0\}$ across $H_\lambda$, while $-\Delta u^r(x,0)=-\partial_y u^r(x,0)\leq 0$ in the reflection of $B^N_R$ across the trace of $H_\lambda$. As for $u$, it still solves the PDE below, and its laplacian is negative when it coincides with $\varphi$. So, this is sufficient to perform the moving plane method until the moving plane touches $\mathrm{supp}~\varphi^+$ . Consider now a point $x_0$ at distance $\e_0>0$ small from $\mathrm{supp}~\varphi^+$, and its projection $\bar x_0$ on $\mathrm{supp}~\varphi^+$ . Set 
 $$e_0=\di\frac{x_0-\bar x_0}{\vert x_0-\bar x_0\vert}.
 $$  Let $H_\lambda$ be given by \eqref{e2.15} with $e=e_0$. By convexity, there is $\lambda_0$ such that $H_{\lambda_0}$ is the supporting hyperplane  of $\mathrm{supp}~\varphi^+$ at $x_0$, and all neighbouring points $x'$ are reached by the reflection across $H_\lambda$, $\lambda-\lambda_0$ small enough. This yields a uniform cone of monotonicity for $u$ directed by $e_0$.
\hfill$\Box$

\begin{rem}
\label{r2.1}
As $\beta'\geq0$, uniqueness holds in \eqref{e2.3}.
\end{rem}

\begin{rem}
\label{r2.3}
For more general $f$'s, another way  to obtain solutions to \eqref{e2.3}  with $u_y\geq0$    is to start from the Cauchy Problem
\begin{equation}
\label{e2.30}
\left\{
\begin{array}{rll}
u_t-\Delta_xu+u_y=&0,\ u(t,x,0)\geq\varphi(x)\quad(t>0,x\in B_R^N)\\
u_t-\Delta u+\beta(u)=&0\quad(t>0,x\in\RR^N,y<0)\\
\ \\
u(t,x,y)=&0\quad(t>0,\vert(x,y)\vert=R,y\leq0)\\
u(0,x,0)=&\varphi(x)\\
u(0,x,y)=&0.
\end{array}
\right.
\end{equation}
The equation on the plane $\{y=0\}$ is a parabolic obstacle problem, however \eqref{e2.30} would yield a solution converging, for large times, to a steady solution. The property $u_y\geq0$ is inherited from the initial datum This would not prove that all solutions $u$ to \eqref{e2.3} satisfy $y_y\geq0$, but would at least have the merit of exhibiting one.
\end{rem}

\noindent  A standard compactness argument now shows, by letting $R\to+\infty$, that a sequence of minimisers of the energy $J_R$ given by \eqref{e2.9} will converge, locally on every compact set, to a solution $u$ of \eqref{e2.1}.  We denote this sequence $(u_{R_n})_n$, we have that the sequence $(J_{R_n}(u_{R_n}))_n$ is nonincreasing, simply because $u_{R_n}$ is an admissible test function for $J_{R_{n+1}}$. This implies that $u$ has finite energy, hence, in particular, $u(X)$ tends to 0 as $\vert X\vert \to\infty$. 

\section{Construction of a solution to \eqref{e1.1}: letting $\e\to0$}
\noindent Before starting the approximation process, let us define an auxiliary solution that will help us in proving that the free boundary on the plane actually exists. Let $\rho_0>0$ be such that 
$B_{\rho_0}^N$ contains $\mathrm{supp}~\varphi^+$, and let $\bar\varphi$ be a spherically symmetric function, equal to $\varphi(0)$ on $B_{\rho_0}^N$, and strictly concave outside $\overline B_{\rho_0}^N$. Call $\bar u_\Lambda(x,y)$ a solution of \eqref{e2.1} with $\varphi$ replaced by $\Lambda\bar\varphi$, with $\Lambda>1$. We claim that $\bar u_\Lambda$ is axially symmetric and, if $\Lambda>1$ is large enough, that $\bar u_\Lambda>u$. This is once again done by the usual sliding type argument, as $u$ may exceed  $\bar u_\Lambda$ only when it is close to 0, that is, at infinity. However, the usual maximum principle applies to $\bar u_\Lambda-u$ in this area. Decreasing back $\Lambda$ to 1, we see that $\bar u_1>u$. This function $\bar u_1$ will be called $\bar u$, until further notice.

\medskip
\noindent We  now set  $\beta(u):=\beta_\e(u)$ given by \eqref{e2.8}, and we call $u_\e$ a so constructed solution to \eqref{e2.1}. We want to pass to the limit $\e\to0$, and the main step is a gradient bound.
\begin{lm}
\label{l2.2}
There is $M>0$ independent of $\e$ such that $\Vert \nabla u_\e\Vert_\infty\leq M$.
\end{lm}
\noindent {\sc Proof.} The situation that we have inherited is the following: for all $\lambda>0$, each level set $\Sigma_\lambda=\{u_\e=\lambda\}$ is a graph in the vertical direction, that meets the plane $\{y=0\}$ because $u_\e(X)$ goes to 0 as $\vert X\vert$ goes to infinity. 

Let us consider a point of $\Sigma_\e$ of the form $(x_\e,-\e)$. Do the classical Lipschitz scaling
$$
v(\xi,\zeta)=\frac{u\biggl((x_\e,-\e)+\e(\xi,\zeta)\biggl)}\e,
$$
the function $v$ solves
 \begin{equation}
 \label{e2.10}
 \left\{
 \begin{array}{rll}
 -\Delta_\xi v+\e v_\zeta=&0\quad (\xi\in\RR^N,\zeta=1)\\
 \ \\
 -\Delta v+\beta(v)=&0\quad(\xi\in\RR^N,\zeta\leq1)\\
 v(0,0)=&1.
 \end{array}
 \right.
 \end{equation}
 We want to prove that $v(0,1)$ is bounded by a universal constant. Let us set $v(0,1)=A$ with $A>1$ and let us examine what this implies. Notice that a crude second derivative bound for $u_\e$ is $\vert D^2u_\e(X)\vert\leq\di\frac{C}\e$, so that, by interpolation,
 $v_\zeta(\xi,\zeta)=v_y(\xi,\zeta)=O(\di\frac1{\sqrt\e})$. Thus the equation for $v$ on the plane $\{\zeta=1\}$ is really $-\Delta v=O(\sqrt\e)$, so that the Harnack inequality yields
 $$
 v(\xi,1)\geq\eta_0A\quad\hbox{for $\vert\xi\vert\leq1$},
 $$
 with $\eta_0>0$ universal. Let $\underline v(\xi,\zeta)$ solve
  \begin{equation}
 \label{e2.11}
 \left\{
 \begin{array}{rll}
  v(\xi,\zeta)=&\eta_0A\quad (\xi\in B_1^N,\zeta=1)\\
  \ \\
 -\Delta v+\Vert f\Vert_\infty v=&0\quad(\xi\in B_1^N,\zeta\leq1)\\
 v(\xi,\zeta)=&0\quad(\xi\in \partial B_1^N,\zeta\leq1)
  \end{array}
 \right.
 \end{equation}
System \eqref{e2.11} implies that $\underline v$ is a subsolution to \eqref{e2.10}, and that $\underline v\leq v$. The argument is of the classical sliding type: we have 
\begin{equation}
\label{e2.12}
\gamma\underline v\leq v
\end{equation}
 for small $\gamma>0$ (one has to examine $v$ and $\underline v$ at infinity, but this is where the standard maximum principle applies) and \eqref{e2.12}  holds on compact sets; this implies that \eqref{e2.12} holds up to $\gamma=1$, if not there would be a touching point between $\gamma\underline v$ and $v$ inside the cylinder $B_1^N\times\{\zeta\leq1\}$. However, \eqref{e2.11} and the strong maximum principle implies that $v(0,0)$ is a universal multiple of $A$, which implies that $A$ is universally bounded.
 By elliptic regularity, $v_\zeta$ and $\nabla_\xi v$ are bounded on $B_1^N\times(-1,1)$.
 
 \noindent Consider $l_0<1$ such that the first eigenvalue of the Dirichlet Laplacian on the cylinder 
 $$\{\xi_1\geq-2,\vert\xi_i\vert\leq 2\}\times(1-2l_0,1)
 $$
  is larger than $1+\Vert\beta_1\Vert_\infty$. Because of the boundedness of $\partial_{\xi_i}v$ and $\partial_\zeta v$ on  $\{\xi_1=-1,\vert\xi_i\vert=2\}\times(1-l_0,1)$, a large multiple of the first eigenfunction of this operator bounds   $\partial_{\xi_i}v$ and $\partial_\zeta v$ in  $\{\xi_1\geq-1,\vert\xi_i\vert\leq1 \}\times(1-l_0,1)$. Because $v$ decreases in the vertical direction, this process actually bounds $v$ in the whole slab $\{\xi_1\geq-1,\vert\xi_i\vert\leq1 \}\times(-\infty,1)$.
 
 \noindent Let $\Sigma_\e$ be the projection of $\Sigma_\e$ onto the plane $\{y=0\}$ and let $\Omega_\e$ the open subset enclosed by $\tilde\Sigma_\e\times\RR_-$.  As the above argument can be repeated at any point where $u=\e$, we may conclude, due to the invariance of the gradient norm under Lipschitz scaling, that $\vert\nabla u_e\vert$ is universally bounded on $\partial\Omega_\e$. 
 
 \noindent Once we are here, it is easy to finish. As $\{u_e=\e\}$ is a graph, it is at distance of the order $\e$ at least from the plane, so that the scaling argument will bound $\nabla u_\e$ on the surface, and below also because $u$ is nonincreasing in the vertical direction. Above the surface, $u_\e$ is harmonic, so that the maximum principle and elliptic regularity yield the  gradient bound. \hfill$\Box$

\noindent  We need a last lemma ensuring that, for the solution $u(x,y)$ to \eqref{e1.1} that we will eventually construct,  the set
 $$
 \mathrm{supp}~u(,.)\backslash\mathcal{C}
 $$
 is nontrivial.
 \begin{lm}
 \label{l3.4}
 With the notations of assumptions 1 and 2 on $\varphi$, there is $\rho_1>0$ and $\delta_1>0$ independent of $\e$ such that, if $\rho_0\leq\rho_1$ we have
 $$
 u_\e(x,0)\geq\delta_1\ \hbox{if}\ \vert x\vert=1+\rho_0.
 $$
  \end{lm}
 \noindent{\sc Proof.}  Pick $\rho_0>0$ such that $\mathrm{supp}~\varphi^+$ is contained in $B_{1+\rho_0/2}^N$. Assume the existence of a sequence  $(\e_n)_n$ going to 0 for which there is a sequence $(x_n)_n$ of $\partial B_{1+\rho_0}^N$ and a sequence $(\delta_n)_n$ going to 0 such that
 $$
 u_{\e_n}(x_n,0)=\delta_n.
 $$
Recall that $u_{\e_n}$ is a minimiser of the energy
 $$
 J_{\e_n}(u)=\int_{\RR^N}\vert\nabla_x u(x,0)\vert^2dx+\int_{\RR^{N+1}_-}\biggl(\vert\nabla u\vert^2+B_{\e_n}(u))\biggl)dxdy,
 $$
 where $B_{\e_n}$ is a primitive of $\beta_{\e_n}$. Let $\Omega_n$ be the $\delta_n$ level surface of $u_{\e_n}$ We have
 $$
  J_{\e_n}(u_{\e_n})\geq \int_{\partial\Omega_n}\vert\nabla_x u_{\e_n}(x,0)\vert^2dx,
  $$
 and the RHS is larger than the infimum, over all functions of $H^1(B^N_{1+\rho_0})$ that are larger than $\varphi$ and equal to $\delta_n$ on $\Omega_n$, in other words, the solution $v_n$ of the classical obstacle problem in $\Omega_n$ with the Dirichlet data $\delta_n$ on $\partial\Omega_n$.  Let us bound its Dirichlet energy from below; we recall that, if $\bar x_n$ is the projection of $x_n$ onto $\mathrm{supp}~\varphi^+$, and
$$
e_n=\frac{x_n-\bar x_n}{\vert x_n-\bar x_n\vert},
$$
then, around $x_n$, the level set $\partial\Omega_n$ is a uniformly Lispchitz graph in the direction $e_n$. In other words, there is a local coordinate system around $\bar x_n$,  named $(\xi',\xi_N)$, a radius $r_0>0$, independent of $\rho_0$, and a function $\psi_n(\xi')$ defined in the ball $B^{N-1}_{r_0}$ of $\RR^{N-1}$, expressed as  $\{\vert\xi'\vert\leq r_0\}$, with $\Vert\nabla\psi_n\Vert_{L^\infty(\{\vert\xi'\vert\leq r_0\})}$ uniformly bounded, such that, in this coordinate system, we have $x_n=(0,\di\frac{\rho_0}2)$, and
$$
\Sigma_n\cap B^{N-1}_{r_0}=\{(\xi',\frac{\rho_0}2+\psi_n(\xi')),\vert\xi'\vert\leq r_0\}.
$$
Let $\Gamma_n$ be the cone with vertex $\bar x_n$, generated by the ball $B^{N-1}_{r_0}$, its aperture  is given by  
$$\theta_n=\mathrm{Arctan}\biggl({2r_0}{\rho_0}\biggl),
$$
thus bounded from below independently of $\rho_0$. And so, there is a universal constant $C_0>0$ such that
$$
\vert\Gamma_n\vert\geq C_0\rho_0.
$$ 
In $\Gamma_n$, the function $v_n$ drops from $v_n=1$ on $B^N_1$ to $O(\e_n)$ on $\Sigma_n$; by the Harnack inequality we have 
$$
-\partial e_n v_n(x)\geq \frac{C}{\rho_0},
$$
$C>0$ universal. 
All in all, we have
 $$
 \int_{\partial\Omega_n}\vert\nabla_x u_{\e_n}(x,0)\vert^2dx\geq\int_{\Gamma_n}\vert\nabla v_n\vert^2dx\gtrsim\frac1{\rho_0},
 $$
 a quantity that goes to infinity as $\rho_0$ goes to 0. On the other hand,  take any function $v(x,y)$ such that
 \begin{enumerate}
 \item $v(.,0)$ is smooth and supported in $B^N_2$,
 \item the set $\partial\{v>0\}$ is a smooth surface that meets the plane $\{y=0\}$ tangentially,
 \item $v$ is supported in $B^{N+1}\cap\{y<0\}$, smooth in $B^{N+1}\cap\{y<0\}$ and decays linearly to 0 in the vicinity of $\partial\{v>0\}$.
 \end{enumerate}
  Then, due to the nondegeneracy of $v$ along  the boundary of its zero set we have, uniformly in $n$:
  $$
   J_{\e_n}(u_{\e_n})=O(1),
 $$
 independently of $\rho_0$. This is a contradiction. \hfill$\Box$

\noindent{\sc Proof of Theorem \ref{t1.1}.} The uniform Lipschitz bound implies the convergence of a subsequence of $(u_\e)_\e$ to  a limiting Lipschitz function $u$. From Lemma \ref{l3.4}, and the fact that $u$ is Lipschitz, the set 
$$
\Omega=\{x\in\RR^N,\ u(x,0)>0\}
$$
is nontrivial, as the maximum of $u$ is $\varphi(0)$. Moreover, if $\mathcal{C}$ is its coincidence set, then $\Omega\backslash\mathcal{C}$ has nonempty interior. This implies that the set
$
\{u>0\}\cap\RR^{N+1}_-
$
has also a nonempty interior. 

\noindent Let $\bar u_\e$ be the axially symmetric barrier constructed for $u_\e$ with the aid of the function $\bar\varphi$. From \cite{CR1} (an analysis based on \cite{AltC} and \cite{BCN}), a subsequence of $(\bar u_\e)_\e$ converges to a minimiser of the functional
$$
\bar J(u)=\int_{\RR^N}\vert\nabla u(x,0)\vert^2dx+\int_{\RR^{N+1}_-}\vert\nabla u\vert^2dxdy+\vert\{u>0\}\cap\RR^{N+1}_-\vert.
$$  
This implies that $\bar u$ develops a free boundary $\bar\Gamma$ that is, due to \cite{AltC}, a $C^\infty$ surface in $\RR^{N+1}_-$ which is, due to axial symmetry, generated by a $C^\infty$ curve $\bar\Gamma_0$. Assume that $\bar\Gamma$ does not meet the plane, integration of the full equation \eqref{e1.1} for $\bar u$ yields
$$
+\infty=\int_{\Gamma}d\sigma(x)=\int_{\mathcal{C}}\partial_y\bar u<+\infty,
$$
a contradiction. This implies that $\bar\Gamma_0$ has finite length, in other words, that $\bar\Gamma$ meets the plane $\RR^N$ along a sphere of finite radius.

\noindent As $\bar u\geq u$, we have 
$$\{u>0\} \subset \{\bar u>0\},
$$
 so that $\Gamma$ lies above $\bar\Gamma$ and intersects the plane $\RR^N$ along $\partial\Omega$. The positivity set $\Omega$ of $u$ in $\RR^N$ is thus bounded. Notice that $\bar u$ was useful in ruling out an infinite surface $\partial\Omega$ with needles at infinity, something that is not {\it a priori} ruled out by the fact that $\{u>0\}$ has finite perimeter.

\noindent We already know that $u(x,0)$ is decreasing in the direction $\di\frac{x}{\vert x\vert}$, as soon as $x$ is outside $\mathrm{supp}~\varphi^+$ . Let us prove that it actually decreases in a cone of directions of axis $\di\frac{x}{\vert x\vert}$. For this, it is enough to prove that the distance between $\partial\Omega$ and $\mathrm{supp}~\varphi^+$  is positive. This, however,  is not difficult to see, as $u(.,0)$ is $C^{1,\alpha}$ for all $\alpha\in(0,1)$. Should there be a coincidence point between $\partial\Omega$ and $\mathrm{supp}~\varphi^+$, this would imply, by the convexity of $\mathrm{supp}~\varphi^+$, a derivative discontinuity at that point. This proves  that the boundary $\partial\Omega$ is a bounded Lipschitz surface  of $\RR^N$ with no boundary.  

\noindent We finally notice that, because $u(x,y)$ is decreasing in a cone of directions of $\RR^{N+1}_-$, $\Gamma\cap\{y<0\}$ is a Lipschitz surface; due to \cite{C2} is is $C^{1,\alpha}$ and the free boundary condition $\vert\nabla u\vert=-u_\nu=1$ ($\nu$ is the exterior normal to $\{u>0\}$) is satisfied. As a consequence of \cite{KNS}, $\Gamma\cap\{y<0\}$ is analytic. Inside $\Omega\backslash\mathcal{C}$, classical ellliptic theory implies that the equation $-\Delta u+u_y=0$ is satisfied in the strong sense in the whole plane $\RR^N$, as well as in the classical sense in $\Omega\backslash\mathcal{C}$.  \hfill$\Box$

\section{Optimal regularity and nondegeneracy}
\noindent Let $u$ be a solution of \eqref{e1.1}. We are now ready to start the analogy with the obstacle problem in the vicinity of $\partial\Omega$, that is, where the free boundary $\Gamma$ generated by the one phase problem in the lower half space $\RR^{N+1}_-$ meets the plane $\RR^N$. To have an idea of what to expect, it is good to recall what happens in the axially symmetric case, that is, when the function $\varphi$ is spherically symmetric in $\RR^N$. In the coordinate system $(r=\vert x\vert, y)$, Problem \eqref{e1.1} reduces to
\begin{equation}
\label{e3.1}
\left\{
\begin{array}{rll}
\Delta u+\di\frac{N-1}ru_r=&0\quad(x,y)\in\{u>0\}\\
\vert\nabla u\vert=&1\quad((x,y)\in\Gamma:=\partial\{u>0\})\\
\ \\
-u_{rr}-\di\frac{N-1}ru_r+ u_y=&0\quad(x\in\RR^N, y=0)\\
u(r,0)\geq&\varphi(r).
\end{array}
\right.
\end{equation}
The positivity set of $u$ on the line $\{y=0\}$ is a segment $[0,R_0)$. In the vicinity of the point $(R_0,0)$ we have \cite{CR1} that $\Gamma$ is a graph in the $r$ variable: 
$$
y=\psi(r),\ r<R_0,\quad\quad\hbox{with}\ \psi(r)=-\frac{(r-R_0)^2}{2}+o_{r\to R_0^-}(R_0-r)^2.
$$
Moreover, we have 
\begin{equation}
\label{e3.2}
\di\lim_{r\to R_0^-}u_y(r,0)=1,
\end{equation} 
so that the ODE for $u(r,0)$, together with the initial data $u(R_0,0)=u_r(R_0,0)=0$ yields
$$
u(r,0)=\frac{(r-R_0)^2}{2}+o_{r\to R_0^-}(R_0-r)^2.
$$
Coming back to the multi-D case, the least one can expect is therefore a quadratic detachment of $\Gamma$ from the plane, as well as a quadratic behaviour of $u$ in the vicinity of $\partial\Omega$. This is what we will
endeavour to prove in this section.

\noindent The first question is now whether $u(x,0)$ enjoys a better regularity than Lipschitz. The obstacle problem suggests $u(.,0)\in C^{1,1}(\RR^N)$, and this is what we will prove. The following theorem is sufficient for our purpose.

\begin{thm}
\label{t3.1}
There is $C>0$ universal such that 
$$
u(x,0)\leq C\vert x\vert^2.
$$
\end{thm}

\noindent{\sc Proof.}  Let $U^x(r)$ be the  average of $u(.,0)$ over the $N$-dimensional sphere of origin $x$ and radius $r$:
$$
U^x(r)=\frac1{r^{N-1}}\int_{\partial B_r^N(x)}u(x',0)d\sigma(x').
$$
We have - notice that this is also valid if $B_r^N(x)$ intersects the coincidence set, although we will not need it:
$$
U^x_{rr}+\frac{N-1}rU^x_r\leq C,\quad U^x(0)=0,
$$
because of the equation $\Delta_{x}u=u_y$ and the boundedness of $u_y$. Because $U^x(r)$ is nonnegative we have also $U^x_r(0)=0$. This entails 
$$
U^x(r)\leq Cr^2,
$$
for a possibly different $C>0$. Translate the picture so that the origin $(0,0)$ is a free boundary point (that is, $(0,0)\in\partial\Omega)$ and $x\in\Omega$. As $u_y\geq0$ and $\varphi$ concave we have
$
\Delta_xu(x,0)\geq0.
$
The mean value formula yields
$$
\begin{array}{rll}
u(x,0)\leq&\di\frac1{\vert x\vert^2}\int_{B^N_{\vert x\vert}(X)}u(x',0)dx'\\
\leq&\di \frac1{\vert x\vert^2}\int_{B^N_{2\vert x\vert}(0)}u(x',0)dx'\\
=&4U^0(\vert x\vert)\\
\leq&C\vert x\vert^2.
\end{array}
$$
This is the theorem. \hfill$\Box$

The next main item to prove is that $u(x,0)$ grows at least quadratically away from a free boundary point. The geometric setting is   the following: $\Gamma$ is a graph over $\Omega$ of the form
\begin{equation}
\label{e3.10}
y=\psi(x),\quad x\in\Omega,
\end{equation}
where $\psi$ is $C^\infty$ in $\Omega$ and $\psi\equiv0$ on $\partial\Omega$, that we recall to be a Lipschitz surface without boundary.

\begin{thm}
\label{t3.2}
Assume the origine to be on the free boundary $\partial\Omega$. There is $C>0$ and $\Lambda>0$ universal such that 
\begin{equation}
\label{e3.20}
u(x,0)\geq C\vert x\vert^2,\quad-\Lambda\vert x\vert^2\leq \psi(x)\leq -\frac1\Lambda\vert x\vert^2.
\end{equation}
if $x\in\Omega$.
\end{thm}
The upper estimate for $\psi$ is straightforward once the estimate for $u(x,0)$ is known, and results from the gradient bound for $u$:
$$
C\vert x\vert^2\leq u(x,0)=u(x,0)-u(x,\psi(x))\leq\Vert u_y\Vert_\infty\vert\psi(x)\vert.
$$
As for the estimate for $u(x,0)$, we will need a better understanding of the behaviour of $u$ and $\Gamma$ near a free boundary point, and the key element is the following equivalent of \eqref{e3.2}.
\begin{thm}
\label{t3.3}
We have
\begin{equation}
\label{e3.4}
\lim_{d(x,\partial\Omega)\to0, x\in\Omega}u_y(x,0)=1.
\end{equation}
\end{thm} 
Theorem \ref{t3.3} implies the lower estimate for $\psi$ in \eqref{e3.20}. It will in turn result from the following
\begin{lm}
\label{l3.1}
Pick $M>0$. Consider a solution $v(x,y)$ in of the free boundary problem in the whole lower half plane
\begin{equation}
\label{e3.5}
\left\{
\begin{array}{rll}
\Delta v=&0\quad(x,y)\in\{u>0\}\\
\vert\nabla v\vert=&1\quad((x,y)\in\Gamma_v:=\partial\{v>0\})\\
\ \\
v(x,0)=&M\quad(x\in\RR^N).
\end{array}
\right.
\end{equation}
Assume the existence of $0<d_1<d_2$ such that 
$$
\Gamma_v\subset\{-d_2\leq y\leq -d_1\}.
$$
Assume the existence of  $\theta_0>0$ such that $\partial_{x/\vert x\vert+e'}v(x,y)\geq0$, for $\vert e'\vert\leq\theta_0$.Then we have 
$$
v(x,y)=(y+M)^+.
$$
\end{lm}
\noindent{\sc Proof.} Consider the comparison function
$$
v_\Lambda(y)=(y+\Lambda)^+,
$$
which is a solution of \eqref{e3.5}. For large $\Lambda>0$, the assumptions on $\Gamma_v$ imply 
\begin{equation}
\label{e3.6}
\{v>0\}\subset\{v_\Lambda>0\}
\end{equation}
Lower $\Lambda$ until one reaches the threshold $\Lambda_0$ where \eqref{e3.6} holds no more. Notice that $\Lambda_0\geq M$> In case of a strict inequality there is a contact point between $v$ and $v_{\Lambda_0}$, either at finite distance or at infinity. The assumptions on $v$ imply that $\Gamma_v$ is Lipschitz, thus smooth. Therefore, in the finite distance case, the standard strong maximum principle and Hopf Lemma preclude this situation, unless $\Lambda_0=M$ and $v\equiv v_{\Lambda_0}$. In the infinite distance case, the usual compactness argument also yields $\Lambda_0=M$ and
\begin{equation}
\label{e3.7}
v(x,y)\leq (y+M)^+.
\end{equation}
Now, for small $\Lambda>0$ we have $v\geq v_\Lambda$. Increase $\Lambda>0$ such that this is not possible anymore: we have in the end $v\equiv v_{M}$ or
\begin{equation}
\label{e3.8}
v(x,y)\geq (y+M)^+.
\end{equation}
Putting \eqref{e3.7} and \eqref{e3.8} together yields the result. \hfill$\Box$

\medskip
\noindent{\sc Proof of Theorem \ref{t3.3}.} It suffices to prove that the detachment of $\Gamma$ from the plane is at most  quadratic, that is, the existence of a universal constant $C>0$ such that, for all $x_0\in\partial\Omega$ we have
\begin{equation}
\label{e3.11}
\vert\psi(x)\vert\leq C\vert x-x_0\vert^2.
\end{equation}
Assume therefore the existence of $x_0\in\partial\Omega$ and a sequence $(x_n)_n$ in $\Omega$  tending to $x_0$ such that 
$$
\lim_{x_n\to x_0}\frac{\vert \psi(x_n)\vert}{\vert x_n-x_0\vert^2}=+\infty.
$$
Set $h_n=\vert \psi(x_n)\vert$ Rescale in a Lipschitz fashion with $h_n$:
$$
v_n(\xi,\zeta)=\frac{u(x_n+h_n\xi,h_n\zeta)}{h_n}.
$$
From Theorem \ref{t3.1}, the sequence $(v_n(.,0))_n$ converges locally uniformly to 0. Below, it develops a free boundary $\Gamma_n$ of the form 
$$
\zeta=\psi_n(\xi),\quad\hbox{with $\psi_n(0)=-1$}.
$$
As $\Gamma_n$ is uniformly Lipschitz, it is uniformly smooth \cite{C2}, so that it converges to a smooth free boundary. This the function $(v_n)_n$ converges, at least in $C^1_{loc}$, to a solution $v_\infty$ of the one-phase free boundary problem in $\RR^{N+1}_-$, while being 0 on the plane $\RR^N$. This implies $\partial_\zeta v_\infty(\xi,0)<0$, a contradiction with the fact that $u_y(x,0)\geq0$. \hfill$\Box$

\noindent {\sc Proof of Theorem \ref{t3.2}.}  
Pick $x_0\in\partial\Omega$, and let $\bar x_0$ be its projection onto $\mathrm{supp}~\varphi^+$; let $e=\di\frac{x_0-\bar x_0}{\vert x_0-\bar x_0\vert}$. We always may assume that $e=e_N$. The solution $u(x,0)$ is decreasing inside a cone directed by $e_N$, so that all its level sets are Lipschitz graphs in the direction $(-e_N)$. Assuming without generality that 
$x_0$ is the origin, it is enough to prove, for $h\in(0,1)$:
\begin{equation}
\label{e3.15}
u(-\sqrt he_N,0)\geq Ch,
\end{equation}
with a universal $C>0$. Rescale in a quadratic way:
$$
U_h(x)=\frac{u(\sqrt{h}e_N,0)}{h},
$$
we already know from optimal regularity  that $(U_h)_h$ is bounded. We would like to show that $\delta_h=U_h(-e_N)$ is bounded from below as $h\to0$.  For all $\delta>0$, let $\Lambda_\delta$ be the $\delta$-level surface of $u$, we have $\Lambda_0=\partial\Omega$. Let $\Lambda_\delta^h$ the $\delta/h$ level surface of $U_h$. The dilation of coordinates preserving the Lipschitz norm, the surfaces $\Lambda_\delta^h$ are globally uniformly Lipschitz, at least as long as $\delta/h$ is less than the supremum of $u$. And so, up to a subsequence, the surfaces $\Lambda_0^h$ and $\Lambda_{\delta_h/h}^h$ converge, respectively, to Lipschitz graphs along the direction $(-e_N)$ that we denote by $\Lambda_0^0$ and $\Lambda^+_0$.  Let $\delta_0$ be the limit of a convergent subsequence of $(\delta_h/h)_h$, it may be 0 or a very small constant, and we have to make sure that it is larger than a universal constant.

\noindent Let $\Omega^0$ be enclosed by the surfaces $\Lambda_0^0$ and $\Lambda^+_0$. Still because of optimal regularity, its width is bounded from below. Up to a subsequence, $(U_h)_h$ converges to a solution $U$   of
$$
\begin{array}{rll}
-\Delta U+1=&0\quad(x\in\Omega^0)\\
U(x)=&0\quad(x\in \Lambda^0_0)\\
 U(x)=&\delta_0\quad(x\in \Lambda^+_0)
\end{array}
$$
If $\delta_0=0$ we have $U<0$ in $\Omega^0$, a contradiction. And so, by continuity, the positivity of $U$ cannot be preserved if $\delta_0$ is below a universal constant. This proves the claim. \hfill$\Box$
\section{H\"older continuity of the vertical derivative up to $\partial\Omega$}
\noindent The task is made easier by the Lipschitz character of $\partial\Omega$, as well as the fact that the free boundary $\Gamma$ below, that we recall to be a graph $\{y=\psi(x)\}$, is a $C^\infty$ surface, until it meets the plane. This will yield an easy upper bound for $u_y(x,0)-1$ as $x$ approaches $\partial\Omega$. A lower bound will be accessible only in an integral sense, as $x$ approaches $\partial\Omega$. Scaling around the points that are well-behaved for $u_y$ will entail the H\"older continuity of $u_y$.
\begin{lm}
\label{l5.1}
There is $\alpha\in(0,1]$ and $C>0$ universal such that, for all $x\in\Omega$ we have
\begin{equation}
\label{e5.1}
u_y(x,0)\leq 1+C\vert x\vert^\alpha.
\end{equation}
\end{lm}
\noindent{\sc Proof.} This lemma involves the specificity of the problem only in a relatively weak manner. Set $v=u_y$, we use once again the Wentzell condition that entails $v_y(x,0)\leq0$ in $\Omega$. Let $\Omega_1$ be a smooth open subset of $\RR^N$ such that 
$$
\mathcal{C}_u\subset\Omega_1\subset\overline\Omega_1\subset\Omega.
$$
For $M>0$ to be large, let $v_{M}$ be the unique solution of the linear problem
\begin{equation}
\label{e5.3}
\left\{
\begin{array}{rll}
-\Delta v_{M}=&0\quad(x\in\Omega\backslash\Omega_1,y\in\RR)\\
v_{M}(x,y)=&1\quad(x\in\partial\Omega,y\in\RR)\\
v_{M}(x,y)=&M\quad(x\in\partial\Omega_1,y\in\RR).
\end{array}
\right.
\end{equation}
We claim that $v_M\geq v$ as soon as we have
\begin{equation}
\label{e5.4}
M>1+\Vert u_y\Vert_{L^\infty(\Omega_1\times\RR_-\cap\{u>0\})}.
\end{equation}
Notice first that $\{u>0\}\subset\Omega\times\RR_-$, so that  it suffices to compare $v$ and $v_M$ on $\Gamma$ and $\Omega\backslash\Omega_1$. As $v_M$ is even in $y$, we have$$
v_y(x,0)\leq\partial_y v_M(x,0)=0\quad\hbox{for $x\in\Omega\backslash\Omega_1$.}
$$
Then, note that $v_M\geq 1$ in its domain of definition, so that 
$$
v(x,y)\leq 1\leq v_M(x,y)\quad\hbox{for $(x,y)\in\Gamma\cap\biggl((\Omega\backslash\Omega_1)\times\RR_-\biggl)$}.
$$
Finally, \eqref{e5.4} implies that $v\leq v_M$ on $\{u>0\}\cap\partial\Omega_1\times\RR_-$. Therefore, the claim $v\leq v_M$ follows from the maximum principle, and so does the lemma: as $\partial\Omega$ is Lipschitz, $v_M$ satisfies \eqref{e5.1} for $C$ and $\alpha$ universal. \hfill $\Box$

\noindent Let us turn to the lower bound. For $x\in\Omega$ close to $\partial\Omega$ let us set $h=\vert x\vert$.
\begin{lm}
\label{l5.2}
There is $\lambda\in(0,1]$, $\alpha\in(0,1]$ and $C>0$ universal such that
\begin{equation}
\label{e5.2}
\frac1{\vert B_{\lambda h}^N\vert}\int_{B_{\lambda h}^N(x)}\vert u_y(x,0)-1\vert dx\leq Ch^\alpha.
\end{equation}
\end{lm}
\noindent {\sc Proof.} Once again because $\partial\Omega$ is Lipschitz we may choose $\lambda>0$ universal such that
$$B_{\lambda h}^N(x)\subset\{x'\in\Omega:\ d(x',\partial\Omega)\geq\frac{h}2\}.
$$
For $q>0$ let $\Sigma_q$ be the piece of cylinder
$$
\Sigma_q=\{u>0\}\cap\biggl(B_k^N(x)\times\RR_-\biggl),
$$
and set
$$
\Gamma_q=\Gamma\cap\Sigma_q.
$$
From Theorem \ref{t3.2}, the set $\Sigma_{\lambda h}(x)$ is included in a cylinder with height of the order no less, and no more than $h^2$. Integrate the equation for $u$ in it. 
We have
$$
\int_{B_{\lambda h}^N(x)}u_y(x',0)dx'=\vert\Gamma_{\lambda h}(x)\vert+O(h^{N+1}).
$$ 
As we have
$$
\vert\Gamma_{\lambda h}(x)\vert\geq\vert B_{\lambda h}^N(x)\vert=\vert B_{\lambda h}^N\vert,
$$
we have
$$
\frac1{\vert B_{\lambda h}^N\vert}\int_{B_{\lambda h}^N(x)}\biggl(u_y(x',0)-1\biggl)dx'\geq Ch,
$$
an estimate that, combined to Lemma \ref{l5.1} and \eqref{l5.1}, implies our lemma. \hfill$\Box$

\medskip
\noindent{\sc Proof of Theorem \ref{t1.2}.} The issue is now to translate Lemmas \ref{l5.1} and \ref{l5.2} into an inequality for the H\"older quotients of $u$ at scales ranging from $h$ to much smaller than $h$. Assume again that $u$ is translated and rotated conveniently so that $0\in\partial\Omega$. Consider $h>0$ small. The preceding lemmas, together with Markov's inequality, imply that
\begin{equation}
\label{e5.5}
\frac{\bigl\vert\{\vert u_y-1\vert\geq h^{\alpha/2}\}\bigl\vert}{\vert B_h^N\vert}\leq Ch^{\alpha/2}.
\end{equation}
We cover $B_h^N$ with balls of radius $h^{1+\alpha/N}$, inequality \eqref{e5.1} guarantees that, if $B_{h^{1+\alpha/N}}^N(x_i)$  is one of these balls, then there is $\tilde x_i\in
B_{h^{1+\frac\alpha{2N}}}^N(x_i)$ such that
\begin{equation}
\label{e5.6}
u_y(\tilde x_i,0)=1+O(h^{\alpha/2}).
\end{equation}
Consider such a $\tilde x_i$ and rescale $u$ around $\tilde x_i$:
\begin{equation}
\label{e5.6}
{u\biggl((\tilde x_i,0)+h^2(\xi,\zeta)\biggl)}={h^2}v_h(\xi,\zeta),
\end{equation}
where we have dropped the index $i$ for commodity, all the more as the upcoming considerations will not depend on $i$. The situation is now the following:
the function $v(\xi,0)$ is bounded, from optimal regularity, in $B_{\alpha/N-1}^N$, while $D^2_\xi v(\xi,0)$ is $O(h^2)$ uniformly in $\xi\in B_{\alpha/N-1}^N$. By interpolation,
$\nabla_\xi v(\xi,0)$ is $O(h)$ uniformly in $\xi\in B_{\alpha/N-1}^N$. From Theorem \ref{t3.3}, we have, for some $M>0$:
$$
v_h(\xi,\zeta)=(M_h+\zeta)^+o_{h\to0}(1).
$$
the constant $M_h$ being universally bounded from above and below. So, we may write
$$
v_h(\xi,\zeta)=\biggl(M_h+\zeta+w_h(\xi,\zeta)\biggl)^+,
$$
the (rescaled) free boundary $\Gamma_h$ being a graph of the form $\{\zeta=-M_h+\varphi_h(\xi)\}$, the function $\varphi_h$ tending to 0 uniformly in $B_{\alpha/N-1}^N$. And so, everything boils down to showing that these $o_{h\to0}(1)$ are actually (small) powers of $h$. A first step towards this is to prove that the quantity
$$
N_h:=\frac{w_h(0,0)}{h^{\frac{\alpha}{2N}}}
$$
is bounded, that is, $v_h(0,\zeta)$ is $h^{\frac{\alpha}{2N}}$-close to $(M_h+\zeta)^+$. Assume this is not true; because $\nabla_xw_h=O(h)$ we have $v_h(\xi,0)\sim h^{\frac{\alpha}{2N}}N_h$ in $B^N_{h^{\alpha/N-1}}$. As $h\to0$ we have, up to a subsequence,
$$
\lim_{h\to0}\frac{w_h}{ h^{\frac{\alpha}{2N}}N_h}=\bar w(\xi,\zeta)
$$
at least in the $C^1$ sense, the function $\bar w$ being harmonic in the strip $\RR^N\times[-M,0]$ ($M$ is a limit of $(M_h)_h$) with $\bar w(\xi,0)\equiv 1$ and $\bar w(\xi,-M)\equiv0$. Thus we have $\bar w_\zeta(0,0)>0$, something that contradicts \eqref{e5.6}. So, $(N_h)_h$ is bounded, so that $w_h$ is uniformly of order $h^{\frac{\alpha}{2N}}$ in $B^N_{h^{\alpha/N-1}}$.

\noindent The second step towards an estimate of $\partial_yu$ is to write $\partial_yw_h$ in integral form. Let $P_h(\xi,\zeta,\xi')$ be the Poisson kernel of the Laplacian in the positivity domain of $u$, remembering that it decays exponentially fast, uniformly in $h$, as $\vert\xi-\xi'\vert$ grows (because the free boundary is bounded and bounded away from the plane), and that its derivative $\partial_yP_h$ has a singularity of the form $\di\frac1{\vert\xi-\xi'\vert^{N+1}}$ uniformly in $h$ we have
$$
\partial_yw_h(\xi,0)\sim\int_{B_{h^{\alpha/N-1}}^N}\partial_yP_h(\xi,0,\xi')(w_h(\xi',0)-w_h(\xi,0)d\xi',
$$
which entails that $\partial_yw_h(\xi,0)$ is of order $h^{\frac\alpha{2N}}$. It is also H\"older continuous, its  H\"older quotient of order $\beta$ ($\beta<1$) being controlled by $h^{\frac\alpha{2N}}$. Scaling back, this proves the  H\"older continuity of $u$ with exponent $\beta<\di\frac\alpha{2N}$. \hfill$\Box$

\noindent {\bf Acknowledgement.} L.A. Caffarelli is supported by NSF grant DMS-1160802. The research of J.-M. Roquejoffre  has received funding from the ERC under the European Union's Seventh Frame work Programme (FP/2007-2013) / ERC Grant Agreement 321186 - ReaDi. He also acknowledges a J.T. Oden fellowship for a long term visit in 2018-19.

\noindent 
{\footnotesize

\end{document}